\input amstex
\documentstyle{amsppt}



\topmatter
\title
{Transfers and Periodic Orbits of Homeomorphisms}
\endtitle
\rightheadtext{Transfers and Periodic Orbits}

\author
Daniel Henry Gottlieb
\endauthor
\leftheadtext{Daniel Henry Gottlieb}

\thanks    
\endthanks

\address
Math. Dept., Purdue University, West Lafayette, Indiana       
\endaddress
\email
gottlieb\@math.purdue.edu      
\endemail
\urladdr
http://www.math.purdue.edu/~gottlieb
\endurladdr
\date 
{\thismonth}{\space\number\day,} {\number\year}
\enddate
\subjclass...\endsubjclass
\abstract...\endabstract
\define\thismonth{\ifcase\month 
  \or January\or February\or March\or April\or May\or June%
  \or July\or August\or September\or October\or November%
  \or December\fi}


\abstract
Bo Ju Jiang applied Neilsen theory to the study of periodic orbits of a homeomorphism. His method
employs a certain loop in the mapping torus of the homeomorphism. Our interest concerns
the persistence of periodic orbits in parameterized  families of homeomorphisms. This leads us to
consider fibre bundles and equivariant maps, which gives us a nice point of view.
\endabstract



\NoBlackBoxes

\endtopmatter
\parskip=6pt
\document 

\head 1.\ Introduction\endhead

Bo Ju Jiang introduced mapping tori into the study of periodic orbits of a homeomorphism
in his article [Jiang(1996)]. If $h: F \rightarrow F$ is a homeomorphism, then a periodic orbit of
length $n$ gives rise to a loop $\sigma$ in the mapping torus $T_h$ which represents 
an element $[\sigma]$ in the fundamental group of the mapping torus $\pi_1(T_h)$. This
element $[\sigma]$ plays a role in the Neilsen theory for the mapping torus and gives information
about the existence of a periodic orbit of the homeomorphism $h$.

Now the mapping torus $T_h$ gives rise to a fibre bundle $p:T_h \rightarrow S^1$ over a circle
with fibre $F$. The loop $\sigma$ can be regarded as an $n$-covering space over the base circle when the fibre bundle projection $p$ is restricted to the image of $\sigma$ in $T_h$.

We are interested in the persistence of periodic orbits under isotopies of $h$, or more generally,
under parameterized families of homeomorphisms. This suggests that we consider fibre bundles
with fibre $F$ over a base $B$ and ask when does the total space $E$ contain a subspace 
$S$ so the the restriction of the fibre bundle projection $p:E \rightarrow B$ to $S$ results in a 
map $S \rightarrow B$ of ```degree'' $n$. 

Note this generalizes the fact that $\sigma$ maps to
$S^1$ with degree $n$. We also point out that when $B$ is not an oriented manifold, we take the
definition of degree of a map given in [Gottlieb(1986)]: The {\it degree} of a map $f:X \rightarrow Y$ is taken
to be the smallest positive integer $N$ for which there exists a homomorphism
$\tau : H_*(Y) \rightarrow H_*(X)$ on integral homology so that the composition $f_* \circ \tau$
is given by multiplication by $N$. If there is no such integer, than the degree is defined to be zero.
The homomorphism $\tau$ will be called a {\it transfer} for $f$ for the purposes of this paper.

In the case of a single fixed point, there is an integer invariant called the index of the fixed point
which guarantees the persistence of a fixed point under homotopy if the index is nonzero. This is
also true for zeros of vectorfields; see [Gottlieb and Samaranayake(1994)]. This persistence
leads to transfers for fibre bundles with vectorfields tangent to the fibre, [Becker and Gottlieb(1991)], and
to transfers for fibre bundles with fibre preserving maps, [Dold(1976)]. 

A more direct transfer related to the Lefschetz number and to Reidemeister torsion was studied by
Wolfgang L\"uck in [L\"uck(1997)], but it involves algebraic $K$-theory groups instead of homology
groups.

We therefore initiate a study of transfers in a setting which parameterizes  Jiang's approach
by recalling relevant facts about principal fibre bundles and their associated principal fibre
bundles. Then we investigate a few simple consequences of the existence of transfers.

\head 2.\ An Example \endhead

Now the circle group $S^1$ acts on itself via a multiplication map
$ \mu : S^1 \times S^1 \rightarrow S^1$ given by 
 $(e^{2 \pi ia} , e^{2 \pi ib}) \mapsto e^{2 \pi i(a+b)}$. Now left multiplication by $e^{2 \pi ia} $ is a 
 self homeomorphism, and it has a periodic orbit of length equal to the denominator of $a$ if $a $
 is a
 fraction with relatively prime numerator and denominator. If $a$ is irrational, there is no periodic
 orbit.
 
 On the other hand, the multiplication map $\mu$ is the clutching map of a fibre bundle over the
 two dimensional sphere $S^2$. The fibre bundle is fibre bundle equivalent to the Hopf  bundle 
 $$
 S^1 \rightarrow S^3 \rightarrow S^2
 $$
 and since the second homology group of the total space is zero and the second homology group
 of the base space is infinite cyclic, the degree of the projection must be zero.
 
 If we alter the map $\mu$ by a homotopy, so that the the map of $S^1$ into the space of 
 self homeomorphisms of $S^1$ corresponding to $\mu$ is homotopied to another map
 $ S^1 \rightarrow Homeo(S^1)$, the corresponding clutching map still is fibre bundle equivalent  to the Hopf bundle, 
 and so the degree is zero and so there is ``no transfer''. Hence we expect that there is no global ``coherent'' periodic orbit for the family of homeomorphisms that the clutching map represents. Later
 in this paper we will describe more precisely what we mean by coherent.
 
 \head 3.\ Principal Bundles and Associated Bundles \endhead
 
 A {\it principal bundle} $G \rightarrow E \rightarrow B$ is a fibre bundle whose fibre is a group $G$ and
 $G$ acts freely on the total space $E$. Hence the quotient map $q: E \rightarrow E/G$ is the
 fibre bundle projection. 
 
 Given an action of $G$ on a space $F$, we consider the diagonal action of $G$ on $E \times F$.
 Then the {\it Borel construction} $E \times_G F$ is defined to be the quotient space of 
 the diagonal action of $G$ on $E \times F$. This results in the {\it associated fibre bundle}
 
 $$F \rightarrow E \times_G F \rightarrow B$$.
 
 We can regard the particular action of $G$ on $F$ as a homomorphism 
 $$\rho : G \rightarrow Homeo(F)$$
 from the group to the space of homeomorphisms of $F$. So a more precise notation for the
 Borel construction is $E \times_\rho F$.
 
 \proclaim{Lemma 1}  A cross-section $s$ to the associated bundle $F \rightarrow E \times_G F \rightarrow B$ corresponds to a $G$-map $\hat s : E \rightarrow F$
\endproclaim

\demo{proof} Let $s: B \rightarrow E \times_\rho F$ be the cross-section. Then
$s(b) = [e, x] = [ge, gx] $ where the projection $p: [e, x]  \mapsto b$. Then we define 
$\hat s : E \rightarrow F$ by $\hat s : e \mapsto x$. This is indeed an equivariant map
since $[e, x] = [ge, gx] $ implies $\hat s (ge) = gx$.

Conversely, given an equivariant map $\hat s$, we define a cross-section by
$s: b \mapsto [e, \hat s(e)]$. This is well defined since $[e, \hat s(e)] = [ge, g \hat s(e)] = [ge, \hat s(ge)]$
\qed
\enddemo

Important examples of this lemma are:

\noindent  (1) If $\rho : G \rightarrow Homeo(G)$ is given by right multiplication on $G$ by itself,
then the associated bundle $ G \rightarrow E \times_\rho G \rightarrow B$ is the principal bundle
itself. Then a $G$-map $E \rightarrow G$ gives a cross-section of the principal bundle. In this case,
the principal bundle is trivial for two different reasons: The fibre is a retract of the total space; and the 
bundle  has a cross-section.

\noindent (2) If $\rho : G \rightarrow Homeo(G)$ is given by conjugation on $G$, then a cross-section
to the associated bundle corresponds to a {\it principal bundle equivalence} $ h: E \rightarrow E$. This is
an equivariant map which induces the identity on the base space $B$. This follows because the
cross-section corresponds to a $G$-map $f:E  \rightarrow G$, where $G$ acts on itself by the adjoint
representation. Then the map $h: e \mapsto f(e)e$ is the corresponding bundle equivalence since
$f(ge) = g f(e)g^{-1}$.

Now Milnor(1956) showed that every topological group $G$ has a {\it universal principal bundle}
$$ G \rightarrow E_G \rightarrow B_G$$ (where $E_G$ is contractible).
This means that the set of principal $G$-bundle equivalence classes over a base space $B$ are in 
one-to one correspondence with the homotopy classes $[B, B_G]$ when $B$ is a CW-complex. The
correspondence is induced by assigning to a map $f : B \rightarrow B_G$ the {\it pullback bundle}
$ G \rightarrow f^*(E_G) \rightarrow B$ where the total space 
$f^*E_G = \{ (b, e) \in B \times E_G \mid f(b) = p(e) \}$.

A related universal fibration holds for Hurewicz fibrations. If 
$$F \rightarrow E_\infty \rightarrow B_\infty$$ is a universal Hurewicz fibration, then let
$E_\infty^{(F)}$ denote the space of maps of $F \rightarrow E_\infty$ which are homotopy equivalences into
fibres of $E_\infty$. Then $E_\infty^{(F)} \rightarrow B_\infty$ is a principal fibration whose fibre is the space of 
self-homotopy-equivalences. Now $E_\infty^{(F)}$ is essentially contractible so $B_\infty$ is the classifying space of the group of self-homotopy-equivalences of $F$. These results have 
topological issues attached to them. The most recent expose is Booth(2000). 

So we may consider the sequence of homomorphisms $G \rightarrow \Cal H \rightarrow \Cal E$ 
where $\Cal H$ denotes the group of homeomorphisms of $F$ and $\Cal E$ denotes the monoid of
self homotopy equivalences of $F$. Then a fibre bundle with fibre $F$ and group $G$ may be
discussed using the sequence of induced maps 
$B  \rightarrow B_G  \rightarrow B_\Cal H \rightarrow B_\Cal E$. That is, given a fibration over
$B$, there is a classifying map $k : B \rightarrow B_\Cal E$. This map factors, up to homotopy,through 
$B_\Cal E$ if and only if the fibration over $B$ is fibre homotopy equivalent to a fibre bundle. Similarly,
if $k$ factors through $B_G$, then the fibre bundle is bundle equivalent to a bundle with structure 
group $G$.

Finally, we can classify the groups of self principal bundle equivalences and fibre homotopy 
equivalences. Let $\Cal G$ denote the topological group of self principal bundle equivalences.
Then [Gottlieb(1972)] showed that for a principal or associated bundle over a base space $B$
whose classifying map is $k: B \rightarrow B_G$, the classifying space of $\Cal G$ is the space
of maps of the base space into the classifying space which are homotopic to $k$. In symbols:
$B_\Cal G = Map(B, B_G ; k)$. The same result hold for self fibre homotopy equivalences,
[Gottlieb(1968)] and [ Gottlieb(1970)]. See also [Booth, Heath, Morgan, Piccinini(1984)]. The group 
$\Cal G$ is also called
the group of {\it gauge transformations} for the relevant principal bundle, [Atiyah, Bott(1983)]. 

\head 4. Mapping Tori \endhead

The universal principal bundle for the integers $\Bbb Z$ is the universal covering space
of the circle, $\Bbb Z \rightarrow \Bbb R \rightarrow S^1$. If $h$ is a self homeomorphism
of $F$, the homomorphism $\rho$ from $\Bbb Z$ to $\Cal H$ given by $n \mapsto h^n$  gives rise to
a map $S^1 = B_\Bbb Z \rightarrow B_\Cal H$. This map pulls back to a principal $\Cal H$
bundle $E$ over the circle. If we apply the Borel construction $E \times _\rho F $ to this bundle,
we get the mapping torus $T_h$. Equivalently, we may take the Borel construction 
$E_\Cal H \times_\Cal H F$ to get the associated bundle with fibre $F$ and then take the pullback by the
classifying map to get the mapping torus.

Now if $G$ acts on $F$ and if a subspace $A$ of $F$ is invariant under the action, we obtain
a sub-bundle $E \times_G A$ contained in $E \times_G F$. In the case of the mapping torus $T_h$,
if $A$ is an orbit of $h$, it is invariant under the action of the representation of $\Bbb Z$ on $F$ and
so we get a one dimensional sub-bundle in $T_h$. This is the origin of the loop $\sigma$ in the mapping torus which plays an important  role in [Jaing(1996)].

Jiang also considers the situation in which there is an isotopy from $h: F \rightarrow F$ to the identity
map. This leads to invariants such as braid groups. The isotopy corresponds to a bundle equivalence
(or a gauge transfomation using different words) $f: T_h \rightarrow T_h$. The fibre bundle
$F \rightarrow T_h \rightarrow S^1$ corresponds to a classifying map
$k : S^1 \rightarrow  B_\Cal H$. The group of gauge transformations $\Cal G$ for this bundle has a
classifying space $B_\Cal G = Map(S^1, B_\Cal H; k)$. The group of gauge transformations 
$\Cal G$ is  homotopy equivalent to the loop space of its classifying space. Hence the path componants
of $\Cal G$ are in one-to-one correspondence to the fundamental group of the classifying space
$\pi_1(B_\Cal G) = \pi_1( Map(S^1, B_\Cal H; k))$. In the case at hand, when $h$ is isotopic to the identity, the map $k$ must be 
homotopic to a constant map and the fibre bundle must be trivial.

Now the unsettled part of Neilsen theory occurs in dimension 2. In the case that $F$ is a closed connected oriented surface which is not the 2-sphere, M. E. Hamstrom has shown that the space
of homeomorphisms has  path components for which the higher homotopy groups are zero.
That is each component is aspherical, [Hamstrom(1965), Hamstrom(1966)].

\proclaim {Lemma 2}
The space of maps $Map(X,Y; k)$ where $Y$ is a $K(\pi,1)$ is itself an aspherical CW complex with fundamental group $\pi_1( Map(X,Y; k))$ equal to the centralizer of the image of  
$k_* : \pi_1(X) \rightarrow \pi_1(Y) = \pi$.
\endproclaim

\demo{proof}
In [Gottlieb(1965)], it is shown that the identity component  $Map(Y,Y; 1_Y)$ is aspherical when
$Y$ is aspherical , and the fundamental group is isomorphic to the center of the fundamental 
group of $Y$. The argument that the mapping space is aspherical does not depend on the 
the domain CW complex. The argument that the fundamental group is the center only depends on the
fact that for any $ \alpha$ in the fundamental group of the domain $X$, its image $k_*(\alpha)$ must commute with all of the fundamental group of $Y$. This is Lemma 2 of [Gottlieb(1969)]. 
 \enddemo
 
 Now suppose $F$ is a $K(\pi,1)$. Then lemma 2 implies that any component of the space
 of self homotopy equivalences is a $K(Z\pi, 1)$ where $Z\pi$ is the center of $\pi$. Hence
 $B_\Cal E$ has fundamental group equal to the group of outer automorphisms $Out(\pi)$, and the second homotopy
 group equal to the center of $\pi$. Thus if $\pi$ has trivial center, $B_\Cal E$ is also aspherical.
Thus the group of self fibre homotopy equivalence classes of the mapping torus of $h$ for F  an aspherical complex, denoted $\pi_0(\Cal G)$, is isomorphic to 
 $\pi_1(Map(S^1, B_\Cal E; k))$ which in turn is isomorphic to the centralizer of the image
 of $k_*: \pi_1(S^1) \cong \Bbb Z \rightarrow \pi_1(B_\Cal E) \cong Out(\pi_1(F))$. 
 
 On the mapping level, this shows that the fibre homotopy classes of fibre homotopy equivalences
 of $T_h$ correspond to those homotopy equivalences $f : F \rightarrow F$ so that
 $f \circ h \sim h \circ f$ where  $\sim$  denotes ``is  homotopic to''.
 
 Similarly, for $F$ a closed connected orientable surface with negative Euler-Poincare number, 
 $B_\Cal H$ is an aspherical space by [Hamstrom(1966)] and the fact that the fundamental group of
 $F$ has trivial center. Hence the fibre isotopy classes of
 gauge transformations of $T_h$ form a group isomorphic to the group of isotopy classes of
 self homeomorphisms of $F$ which commute with $h$.
 
 Where as an orbit of $h$ in $F$ gives rise to a subbundle of $T_h$, a gauge transformation
 does not preserve that subbundle in general, but embeds it into a homeomorphic subbundle.
 This situation must give rise to embedding type invariants, such as braid groups.
 
 \head 5. Bundles over Spheres \endhead
 
 The obvious generalization of mapping tori are bundles over spheres. Here the homeomorphism
 $h:F \rightarrow F$ is replaced by the clutching map $c: F \times S^{n-1} \rightarrow F $ (which corresponds to
 $\hat c : S^{n-1} \rightarrow \Cal H$), and the unit interval is
 replaced by the unit $n$-ball $D^n$. Then identifying $F \times D^n$ to $E'$ by using the clutching
 map on $F \times \partial D^n$, we obtain a fiber bundle $F \rightarrow E' \rightarrow S^n$. Now 
 this fibre bundle has an associated principal bundle $ \Cal H \rightarrow E \rightarrow S^n$. So the
 associated bundle $E \times_ \Cal H F$ is fibre bundle equivalent to $E'$.
 
 Now if $\hat c$ maps into a subgroup $G$ of $\Cal H$ which preserves a subspace $A \subset F$, then 
 $E \times_G A$ is a subbundle of  $E'$. If $A$ is a point, then the fibre of this subbundle is a point,
 and so there is a cross-section to $E' \rightarrow S^n$. Thus by lemma 2, there is a $G$-map 
 from $E$ to $F$. 
 
 Not all cross-sections of $E'$ correspond to a $G$-subbundle. But every one corresponds to a
 $G$-map  $E \rightarrow F$.
 
 Now if $E'$ contains an oriented closed submanifold $M$ of dimension $n$ which maps onto $S^n$ 
 under the projection with degree $m$, then the pullback $p^*:E'$  of $E' \rightarrow S^n$ over $M$ has
 a cross-section. So $p^*E'$ has a $G$ map into $F$.
 
 These submanifolds give rise to transfers of $E' \rightarrow S^n$ of trace $m$. In case we have
 a fibre preserving map $f:E' \rightarrow E'$ inducing the identity on the base with Lefschetz number
 $\Lambda_f$, there will be such a manifold, [Gottlieb(1977)]. 
 
 A corollary of this is: For every fibre bundle over a base space which is a closed oriented
 manifold (such as $S^n$) and whose fibre is a closed manifold whose fibre has Euler-Poincare
 number equal to some nonzero  integer $m$ , then there is a manifold which maps onto the
 base space by a map of degree $m$ so that the pullback bundle has a cross-section. Hence
 the associated principal bundle's total space maps equivariantly into the fibre $F$.


\enddocument